\documentclass[12pt]{article}
\usepackage{amssymb}
\usepackage{amscd}
\usepackage{amsmath}
\usepackage{amsfonts}

\newcommand{\eq}{\begin{equation}}
\newcommand{\en}{\end{equation}}
\newtheorem{proposition}{Proposition}
\newtheorem{theorem}{Theorem}
\newtheorem{lemma}{Lemma}

\def\BP{{\Bbb P}}
\def\BP{{\Bbb P}}

\def\BE{{\Bbb E}}

\def\calu{{\cal U}}

\def\calm{{\cal M}}
\def\caln{{\cal N}}

\def\proof{\noindent{\bf Proof\ \ }}
\def\qed{\mbox{\rule{0.5em}{0.5em}}}

\def\cals{{\cal S}}
\def\bfa{{\bf a}}

\newcommand{\ignore}[1]{}

\begin{document}
%\baselineskip=24truept

\title{Asymptotic enumeration of 2-covers and line graphs}
\author{Peter Cameron, Thomas Prellberg and Dudley Stark\\
\small{School of Mathematical Sciences}\\
\small{Queen Mary, University of London}\\
\small{Mile End Road, London, E1 4NS  U.K.}}
\date{}

\maketitle

\begin{abstract}
In this paper we find asymptotic enumerations for the number of line graphs on $n$-labelled
vertices and for different types of related combinatorial objects called 2-covers.

We find
that the number of 2-covers, $s_n$, and proper 2-covers, $t_n$, on $[n]$ both
have asymptotic growth
$$
s_n\sim t_n\sim B_{2n}2^{-n}\exp\left(-\frac12\log(2n/\log n)\right)= B_{2n}2^{-n}\sqrt{\frac{\log
n}{2n}},
$$
where $B_{2n}$ is the $2n$th Bell number,
while the number of restricted 2-covers, $u_n$,
restricted, proper 2-covers on $[n]$, $v_n$, and line graphs $l_n$,
all have growth
$$
u_n\sim v_n\sim l_n\sim B_{2n}2^{-n}n^{-1/2}\exp\left(-\left[\frac12\log(2n/\log
n)\right]^2\right).
$$

In our proofs we use probabilistic arguments for
the unrestricted types of 2-covers and
 and generating function methods for the restricted types of 2-covers
and line graphs.
\end{abstract}

\noindent{\sc keywords: asymptotic enumeration, line graphs, set partitions}

\section{Introduction}

A $k$-cover of $[n]:=\{1,2,\ldots n\}$ is a multiset
of subsets $\{S_1,S_2,\ldots, S_m\}$, $S_i\subseteq [n]$,
(possibly with $S_i=S_j$ for some $i\neq j$),
such that for each $d\in [n]$ the
number of $j$ such that $d\in S_j$ is exactly $k$.
A $k$-cover is called {\em proper} if $S_i\neq S_j$ whenever $i\neq j$.
A $k$-cover is called
{\em restricted} if the intersection of any $k$ of the $S_i$ contains
at most one element. These definitions have been taken from \cite{GJ}.
Note that for a proper $k$-cover $\{S_1,\ldots,S_m\}$ is a set.

The \emph{line graph} $L(G)$ of a simple graph $G$ is the graph whose vertex set is the edge set
of $G$ and such that two vertices are adjacent in $L(G)$ if and only if the corresponding edges of
$G$ are adjacent.

Let $s_n$ be the number of $2$-covers of $[n]$;
let $t_n$ be the number of proper 2-covers of $[n]$;
 let
$u_n$ be the number of restricted, proper 2-covers of $[n]$;
let
$v_n$ be the number of restricted, proper 2-covers of $[n]$;
and let $l_n$ be the number of line graphs on $n$ labelled vertices.
Let $B_n$ be the $n$th Bell number.
Given sequences $a_n$ and $b_n$, we write $a_n\sim b_n$ to mean
$\lim_{n\to\infty} a_n/b_n=1$.
\begin{theorem}\label{main}
The number of $2$-covers and the number of proper $2$-covers have asymptotic growth
\eq\label{tasymp} s_n\sim t_n\sim B_{2n}2^{-n}\exp\left(-\frac12\log(2n/\log n)\right) \en while
the number of restricted 2-covers, restricted, proper 2-covers and line graphs all have asymptotic
growth \eq\label{uasymp} u_n\sim v_n\sim l_n\sim
B_{2n}2^{-n}n^{-1/2}\exp\left(-\left[\frac12\log(2n/\log n)\right]^2\right). \en
\end{theorem}

We make some initial observations regarding 2-covers, special graphs and orbits in Section~2. We
use a probabilistic method to prove (\ref{tasymp}) in Section 3. A pair of technical lemmas are
proven in Section~3.1, (\ref{tasymp}) is proven for $s_n$ in Section~3.2 and it is proven for
$t_n$ in Section~3.3. We prove (\ref{uasymp}) in Section~4.

In both probabilistic and generating function proofs we will make use of
Lambert's $W$-function $W(t)$,
 which is a solution to
\eq\label{Lambertdef}
W(t)e^{W(t)}=t
\en
and which has asymptotics (see (3.10) of \cite{MW})
\eq\label{Lambertasymp}
W(t)=\log t-\log\log t+\frac{\log\log t}{\log t}+o\left(\frac1{\log t}\right)\quad\mbox{as}\quad t\rightarrow\infty\,.
\en

For each $k$-cover $S_1,\ldots,S_m$ of $[n]$
we can define an associated $m\times n$
incidence matrix $M$
with entries given
by
$$
M_{i,j}=
\left\{
\begin{array}{l l}
1&{\rm if \ }j\in S_i;\\
0&{\rm if \ }j\not\in S_i.
\end{array}
\right.
$$
Note that $M$ has exactly $k$ ones in each column and that the rows are unordered. A $k$-cover is
proper if and only if $M$ has no repeated rows. A $k$-cover is restricted if and only if $M$ has
no repeated columns. Therefore, Theorem~\ref{main} is equivalent to the asymptotic enumeration of
certain 0-1 matrices. The general methods of this paper were used for the asyptotic enumeration of
other 0-1 matrices called incidence matrices in \cite{CPS1,CPS2}.

\section{$2$-covers, line graphs and orbits}
In this section we establish correspondences between $2$-covers, line graphs and orbits of certain
permutation groups.

\subsection{$2$-covers and graphs}
We define a {\em special multigraph} to be a multigraph with no isolated vertices or loops. Our
first result is
\begin{proposition}\label{specialm}
There is a bijection between $2$-covers on $[n]$ and special multigraphs having unlabelled
vertices and  $n$ labelled edges, such that
\begin{itemize}
\item proper $2$-covers correspond to multigraphs having no connected
component of size~$2$;
\item restricted $2$-covers correspond to simple graphs.
\end{itemize}
\end{proposition}
\proof
Let $S_1,\ldots,S_m$ be a $2$-cover of $[n]$. Construct a graph $G$ as
follows:
\begin{itemize}
\item the vertex set is $[m]$;
\item for each $i\in[n]$, there is an edge $e_i$ joining vertices $j$ and
$k$, where $S_j$ and $S_k$ are the two sets of the $2$-cover containing $i$.
\end{itemize}
The graph $G$ is a multigraph (that is, repeated edges are permitted), but it
has no isolated vertices and no loops.

Conversely, given a multigraph without isolated vertices or loops, we can recover a $2$-cover:
number the edges $e_1,\ldots,e_n$, and let $S_i$ be the set of indices $j$ for which the $i$th
vertex lies on edge $e_j$. Thus we have the first part of the proposition.

The second part comes from observing that a ``repeated set'' in a $2$-cover corresponds to a pair
of vertices lying on the same edges, while a pair of elements lying in two different sets
correspond to a pair of edges incident to the same two vertices.
\hfill\qed

\subsection{Generating function identities for $2$-covers}

Recall that $s_n$, $t_n$, $u_n$ and $v_n$ denote the numbers of $2$-covers, proper $2$-covers,
restricted $2$-covers, and restricted proper $2$-covers respectively. Using
Proposition~\ref{specialm} in this subsection we will find relationships between these quantities
and derive corresponding generating function identities.

\begin{proposition}\label{ident}
Let $S(n,k)$ denote the Stirling numbers of the second kind, that is,
the number of set partitions of $[n]$ into exactly $k$ nonempty subsets.
Then,
\begin{eqnarray*}
s_n &=& \sum_{k=1}^n S(n,k) u_k \\
t_n &=& \sum_{k=1}^n S(n,k) v_k \\
u_n &=& \sum_{k=0}^n {n\choose k} v_k
\end{eqnarray*}
\end{proposition}

\proof
We prove these for the corresponding special multigraphs.

Any special multigraph with edges $e_1,\ldots, e_n$ can be described by
giving a partition of $[n]$ into, say, $k$ parts, together with a special
simple graph with $k$ labelled edges; simply replace the $i$th edge of the
simple graph by the $i$th set of edges of the partition (where
the edges are ordered lexicographically, say). This is clearly a
bijection. Moreover, the simple graph has no connected components of size~$2$
if and only if the same holds for the multigraph. This proves the first two
equations.

Given a special simple graph, there is a distinguished subset of $[n]$ (of size $n-k$, say)
consisting of isolated edges; the remaining graph has no components of size~$2$. Again, the
correspondence is bijective. So the third equation holds.
\hfill\qed

Proposition~\ref{ident} can be reformulated in terms of exponential generating functions. Let
$S(x)=\sum_{n\ge0}s_nx^n/n!$, with similar definitions for the others. The proof of
Proposition~\ref{genfuncs} is omitted.
\begin{proposition}\label{genfuncs}
\begin{eqnarray*}
S(x) &=& U(e^x-1) \\
T(x) &=& V(e^x-1) \\
U(x) &=& V(x)e^x.
\end{eqnarray*}
\end{proposition}

It follows from Proposition~\ref{genfuncs} that $S(x)=T(x)B(x)$, where $B(x)=e^{e^x-1}$ is the
exponential generating function for the Bell numbers. This is easily proved directly.

\subsection{Unrestricted $2$-covers and orbits}

Recall the notation $F_n(G)$ for the number of orbits of the
oligomorphic group $G$ on ordered $n$-tuples of distinct elements,
and $F_n^*(G)$ for the number of orbits on all $n$-tuples. Let
$S_\infty^{\{2\}}$ denote the group induced by the infinite
symmetric group on the set of all $2$-element subsets of its domain.

\begin{proposition}
$F_n(S_\infty^{\{2\}})=u_n$ and $F_n^*(S_\infty^{\{2\}})=s_n$.
\end{proposition}
\proof Simply observe that an $n$-tuple of distinct $2$-sets is the edge set of a special simple
graph with $n$ labelled edges, while an arbitrary $n$-tuple of $2$-sets is the edge set of a
special multigraph with $n$ labelled edges. \hfill\qed

We note that the relation
\[F^*(G)=\sum_{k=1}^nS(n,k)F_k(G)\]
gives an alternative proof of the first equation in Proposition~\ref{ident}. We do not know of a
similar interpretation of the other two parameters.

\subsection{Generating function identities for line graphs}
Let $L(x)=\sum_{n\geq 0}l_nx^n/n!$. We now prove
\begin{proposition}\label{line}
\[L(x) = e^{-x^3/3!}U(x) = e^{x-x^3/3!}V(x).\]
\end{proposition}
\proof According to Whitney's Theorem \cite{H}, an isomorphism between line graphs $L(G_1)$ and
$L(G_2)$ of connected graphs is induced by an isomorphism from $G_1$ to $G_2$, except in one case:
the line graphs of the triangle $K_3$ and the star $K_{1,3}$ are isomorphic.

Now the connected components of line graphs which are triangles contribute a factor $e^{x^3/3!}$
to the exponential generating function $L(x)$ for line graphs on $[n]$; that is,
$L(x)=e^{x^3/3!}W'(x)$, where $W'(x)$ is the e.g.f. for line graphs with no such components.
Similarly, components which are triangles or stars contribute a factor $(e^{x^3/3!})^2$ to the
e.g.f. for special simple graphs with $n$ edges. Proposition \ref{line} now follows by Whitney's
Theorem and Proposition~\ref{genfuncs}. \hfill\qed

\section{Unrestricted 2-covers: a probabilistic approach}
In this section we prove (\ref{tasymp}) of Theorem~\ref{main} by using a probabilistic
construction.
\subsection{Technical results}
We proceed with the following definitions and lemma.
Let $T_n$ be the set of proper
2-covers on $[n]$. Let $\cals_n$ be the set of set partitions of $[2n]$.
Let $E_{1,n}\subset\cals_n$ be the subset of set partitions of $[2n]$ such
that $j$ and $j+n$ are contained in different blocks for each $j\in [n]$.
Define the function $\psi$ from a subset $\tilde{S}$ of $[2n]$ to a
subset of $[n]$ by $\psi(\tilde{S})=\{j:j\in\tilde{S}{\rm \ or \ }
j+n\in\tilde{S}\}$. Let $E_{2,n}\subset\cals_n$
be the subset of set partitions of $[2n]$
with blocks
$\{\tilde{S_1},\ldots,\tilde{S_m}\}$
such
that $\psi(\tilde{S}_{i_1})\neq \psi(\tilde{S}_{i_2})$ for each
$i_1\neq i_2$.
Let $C_n=E_{1,n}\cap E_{2,n}$. Let $\phi$ be the function on $\cals_n$
given by
$$
\phi(\{\tilde{S_1},\ldots,\tilde{S_m}\})
=
\{\psi(\tilde{S_1}),\ldots,\psi(\tilde{S_m})\}.
$$
\begin{lemma}\label{proper}
$\phi$ maps $C_n$ onto $T_n$ and $|\phi^{-1}(\bfa)|=2^n$
for all $\bfa\in T_n$.
\end{lemma}
\proof Fix $\{\tilde{S_1},\ldots,\tilde{S_m}\}\in C_n$.
Each $j\in [n]$ appears in exactly two blocks of
$\phi(\{\tilde{S_1},\ldots,\tilde{S_m}\})$ because of the definition
of $E_{1,n}$ and the blocks of $\{\tilde{S_1},\ldots,\tilde{S_m}\}$
are unique because of the definition of
$E_{2,n}$ so $\phi(\{\tilde{S_1},\ldots,\tilde{S_m}\})\in T_n$.

Let $\bfa=\{S_1,\ldots,S_m\}\in T_n$.
For each $j\in [n]$ there are
two ways of assigning $j$, $j+n$ to the appearances of $j$ in $\bfa$
(think of a fixed ordering of the blocks of $\bfa$ to see this).
The choices made for every $j\in [n]$ determine an {\em assignment}.
Clearly, every element of $\phi^{-1}(\bfa)$ must be of the form
$\chi(\bfa)$ for some assignment $\chi$.
There are $2^n$ assignments. We also write $\chi(S_i)$
for the block $\tilde{S_i}$ corresponding to $S_i$ in $\chi(\bfa)$.

We claim that each assignment $\chi(\bfa)$ gives a unique element of $C_n$. To see this, first
note that $j$ and $j+n$ are clearly in different blocks of $\chi(\bfa)$, so $\chi(\bfa)\in
E_{1,n}$. Secondly, $\phi\circ\chi$ is the identity map on $T_n$. Therefore, $\chi(\bfa)\in
E_{2,n}$ because $\bfa$ is a proper 2-cover. Moreover, $\chi_1(\bfa_1)\neq\chi_2(\bfa_2)$ for all
$\bfa_1,\bfa_2\in T_n$ such that $\bfa_1\neq\bfa_2$ and for all assignments $\chi_1$ and $\chi_2$,
which gives $\phi^{-1}(\bfa_1)\cap\phi^{-1}(\bfa_2)=\emptyset$.

We next prove that if $\chi_1$ and $\chi_2$ are two assignments such that
$\chi_1(\bfa)=\chi_2(\bfa)$, then $\chi_1=\chi_2$. To see this, let
$$
\calu=\{j\in [n]: \chi_1{\rm \ and \ }\chi_2 {\rm \ differ \ for \ }
j\}.
$$
Without loss of generality, assume that $j\in S_1$ and $j\in S_2$. Then, either $j\in\chi_1(S_1)$
and $j\in\chi_2(S_2)$ or $j+n\in\chi_1(S_1)$ and $j+n\in\chi_2(S_2)$ It follows that
$\chi_1(S_1)=\chi_2(S_2)$. Therefore, $\phi\circ\chi_1(S_1)=\phi\circ\chi_2(S_2)$ or $S_1=S_2$
violating the assumption that $\bfa$ is proper. We conclude that $\calu=\emptyset$ and that
$\chi_1=\chi_2$. This implies that $|\phi^{-1}(\bfa)|=2^n$. \hfill\qed

Next we generalize Lemma~\ref{proper} to (possibly) improper covers. Let $U_n$ denote the set of
$2$-covers of $[n]$.
\begin{lemma}\label{improper}
$\phi$ maps $E_{1,n}$ onto $U_n$.
Let $\bfa=\{S_1,S_2,\ldots,S_m\}$ be a 2-cover of $[n]$. Let
$\calm$ be the set of $i\in [m]$ such that there does not exist any $j\in [m]\setminus\{i\}$,
$S_j=S_i$. Let
$$
\rho=\frac{m-|\calm|}{2}
$$
be the
number of pairs $\{i,j\}$ such that $S_i=S_j$. Then
$$
|\phi^{-1}(\bfa)|=2^{n-\rho}.
$$
\end{lemma}
\proof Clearly $\phi$ maps $E_{1,n}$ onto $U_n$. Let $\caln=[n]\setminus\{\cup_{i\in\calm}S_i\}$.
Then $\{S_i:i\in\calm\}$ is a proper cover of $\caln$ and Lemma~\ref{proper} implies that
$$
|\phi^{-1}(\{S_i:i\in\caln\})|=2^{|\caln|}.
$$
For each pair $S_{i_1}$, $S_{i_2}$
such that $i_1\neq i_2$ and $S_{i_1}=S_{i_2}$, it must
be true that $\phi^{-1}(S_i)$ consists of two sets $\tilde{S}_1$
and $\tilde{S}_2$ such that for each $j\in S_{i_1}$ either
$j\in \tilde{S}_{i_1}$ and $j+n\in \tilde{S}_{i_2}$ or
$j+n\in \tilde{S}_{i_1}$ and $j\in \tilde{S}_{i_2}$.
The number of choosing unordered sets
$\tilde{S}_{i_1}$, $\tilde{S}_{i_2}$
is $2^{|S_{i_1}|-1}$. Therefore,
$$
|\phi^{-1}(\bfa)|=
2^{|\caln|}\prod 2^{|S_{i_1}|-1}
=
2^{n-\rho},
$$
where the product is over pairs $i_1,i_2$ such that $i_1\neq i_2$
and $S_{i_1}=S_{i_2}$.
\hfill\qed

\subsection{Asymptotic enumeration of proper 2-covers}

From Lemma~\ref{proper} we conclude that $|C_n|=2^nt_n$ so \eq\label{tnexp}
t_n=2^{-n}|C_n|=2^{-n}\frac{|C_n|}{B_{2n}}B_{2n} \en where $B_{2n}$ is the $2n$th Bell number.

We will now prove
\begin{lemma}\label{asymp}
\eq\label{e1} \frac{|E_{1,n}|}{B_{2n}}\sim \sqrt{\frac{\log n}{2n}} \en and \eq\label{e2}
\frac{|E_{2,n}|}{B_{2n}}=1-O\left(\frac{\log^2 n}{n}\right). \en
\end{lemma}
\proof To prove (\ref{e1}), choose an element of $\cals_n$ uniformly at random and let $X$ be the
number of $j\in [n]$ for which $j$ and $j+n$ are in the same block. We have
\eq\label{step1}
\BP(X=0)=\frac{|E_{1,n}|}{B_{2n}}.
\en
We have $X=\sum_{j=1}^n I_j$ where $I_j$ is the indicator random variable
that $j$ and $j+n$ are in the same block.
The $r$th falling moment of $X_n$ is
\begin{eqnarray*}
\BE(X)_r&=&\BE X(X-1)\cdots(X-r+1)\\
&=& \sum \BE(I_{j_1}I_{j_2}\cdots I_{j_r})
\end{eqnarray*}
where the sum is over $(j_1,\ldots,j_r)$ with no repetitions. To find $\BE(I_{j_1}I_{j_2}\cdots
I_{j_r})$ we take $[2n]\setminus\{j_1,j_2,\ldots,j_r\}$ and form a set partition. We then add
$j_k$ to the block containing $j_k+n$ for each $k\in [r]$. This process is uniquely reversible.
Therefore,
$$
\BE(X)_r=\frac{(n)_r B_{2n-r}}{B_{2n}}.
$$

We apply the formula in Corollary 13, page 18, of \cite{Bo} to obtain \eq\label{step2}
\BP(X=0)=\sum_{r=0}^\infty (-1)^r \frac{\BE(X)_r}{r!}=\sum_{r=0}^\infty
\frac{(-1)^r}{r!}\frac{(n)_r B_{2n-r}}{B_{2n}}. \en To analyze (\ref{step2}) we use the expansion
of the Bell numbers \cite{MW,TP}
\begin{eqnarray*}
\log B_n&=&
e^w(w^2-w+1)-{1\over2}\log(1+w)-1-{w(2w^2+7w+10)\over24(1+w)^3}e^{-w}\\
&&-{w(2w^4+12w^3+29w^2+40w+36)\over48(1+w)^6}e^{-2w}+O(e^{-3w})\nonumber\;,
%\\&+&{w(16w^8+164w^7+616w^6+818w^5-1096w^4
%-5793w^3-10032w^2-12208w-12048)\over5660(1+w)^9e^{3w}}\nonumber\\
%&+&O(e^{-4w})\nonumber
\end{eqnarray*}
where $w=W(n)$ is given by (\ref{Lambertdef}), (\ref{Lambertasymp}),
from which we obtain (using Maple)
$$
\log B_{n-r} - \log B_n =
-rw + \frac{rw}{2n}\left(\frac{r}{w+1}+\frac{1}{(w+1)^2}\right)
+ O\left(\frac{r^3w}{n^2}\right).
$$
In particular,
$$
\frac{B_{n-1}}{B_n}\sim \frac{\log n}{n}
$$
so there exists a constant $C>0$ such that
\eq\label{ratiobound}
\frac{B_{n-r}}{B_n}\leq \frac{(C\log n)^r}{(n)_r}.
\en
Moreover,
\begin{eqnarray*}
\log B_{2n-r} - \log B_{2n} &=&
-rv + \frac{rv}{4n}\left(\frac{r}{v+1}+\frac{1}{(v+1)^2}\right)
+ O\left(\frac{r^3v}{n^2}\right)\\
&=&
-r\log n + rc_n+r^2d_n +
+ O\left(\frac{r^3\log n}{n^2}\right),
\end{eqnarray*}
where $v=W(2n)$ has the expansion
$$
v= \log n - \log\log n + \log 2 +\frac{\log\log n}{\log n} -
\frac{\log 2}{\log n} +
o\left(\frac{1}{\log n}\right),
$$
where
\begin{eqnarray*}
c_n&=&\log n - v -\frac{rv}{4n(v+1)^2} \\
&=& \log \log n - \log 2 - \frac{\log\log n}{\log n} +
\frac{\log 2}{\log n} +
o\left(\frac{1}{\log n}\right)
\end{eqnarray*}
and where
$$
d_n=O\left({1\over n}\right).
$$

Using (\ref{ratiobound}) we estimate
\begin{eqnarray}\label{step3}
\left|\sum_{r>\log^{3/2} n}(-1)^r \frac{\BE(X)_r}{r!}\right|
&\leq&
\sum_{r>\log^{3/2} n}\frac{(n)_rB_{2n-r}}{r!B_n}\nonumber\\
&\leq&
\sum_{r>\log^{3/2} n}\frac{(C\log 2n)^r}{r!}\nonumber\\
&=&
(2n)^C\sum_{r>\log^{3/2} n}e^{-C\log 2n}\frac{(C\log 2n)^r}{r!}\nonumber\\
&=&o(1).
\end{eqnarray}

For $r\leq \log^{3/2} n$, we have
$$
\frac{B_{n-r}}{B_n}
=n^{-r}
\exp\left(
rc_n+r^2d_n
+ O\left(\frac{\log^9 n}{n^2}\right)\right)
$$
and
$$
(n)_r=n^r\exp\left(O\left({r^2\over n}\right)\right),
$$
hence
$$
\BE(X)_r=\exp\left(
rc_n+r^2d_n
+ O\left(\frac{\log^9 n}{n^2}\right)\right).
$$
Therefore,
\begin{eqnarray}\label{step4}
\sum_{0\leq r\leq\log^{3/2} n}(-1)^r \frac{\BE(X)_r}{r!}
&=&
\sum_{0\leq r\leq\log^{3/2} n}\frac{(-1)^r}{r!}e^{rc_n+r^2d_n}
\left(1+O\left(\frac{\log^9 n}{n^2}\right)\right)\nonumber\\
&=& \sum_{0\leq r\leq\log^{3/2} n}\frac{(-1)^r}{r!}e^{rc_n} \left(1+d_n r^2 + O\left(\frac{\log^9
n}{n^2}\right) \right)\nonumber\\
&=&\sum_{0\leq r\leq\log^{3/2} n}\frac{(-1)^r}{r!}e^{rc_n}+d_n\sum_{0\leq r\leq\log^{3/2}
n}\frac{(-1)^rr^2}{r!}e^{rc_n}\nonumber\\
&&\,+\left(\frac{\log^9 n}{n^2}\right)\sum_{0\leq r\leq\log^{3/2} n}\frac{e^{rc_n}}{r!}.
\end{eqnarray}

We proceed to approximate the terms in (\ref{step4}). First, we find that
\begin{eqnarray}\label{step5}
\sum_{0\leq r\leq\log^{3/2} n}\frac{(-1)^r}{r!}e^{rc_n}
&=&
\exp\left(-e^{c_n}\right) +
O\left(\sum_{\log^{3/2} n\leq r\leq n}\frac{e^{rc_n}}{r!}\right)\nonumber\\
&=&
\exp\left(-\frac{\log n}{2}\left[1-\frac{\log\log n}{\log n} +
\frac{\log 2}{\log n} +
o\left(\frac{1}{\log n}\right)
\right]\right)
+o(n^{-1/2})\nonumber\\
&\sim&
\sqrt{\frac{\log n}{2n}}.
\end{eqnarray}
We estimate
\begin{eqnarray}\label{step6}
&& d_n\left| \sum_{0\leq r\leq\log^{3/2} n}\frac{(-1)^r}{r!}r^2e^{rc_n}
\right|\nonumber\\
&=&
d_n\left|
\sum_{2\leq r\leq\log^{3/2} n}\frac{(-1)^r}{(r-2)!}e^{rc_n}
+
\sum_{1\leq r\leq\log^{3/2} n}\frac{(-1)^r}{(r-1)!}e^{rc_n}
\right|\nonumber\\
&=&
d_n\left|
e^{2c_n}\sum_{2\leq r\leq\log^{3/2} n}\frac{(-1)^r}{(r-2)!}e^{(r-2)c_n}
+
e^{c_n}\sum_{1\leq r\leq\log^{3/2} n}\frac{(-1)^r}{(r-1)!}e^{(r-1)c_n}
\right|\nonumber\\
&=& d_n\left( \exp\left(-e^{c_n}+2c_n\right) + \exp\left(-e^{c_n}+c_n\right) +
O\left(e^{2c_n}\sum_{\log^{3/2} n\leq r\leq n}\frac{e^{rc_n}}{r!}\right)
\right)\nonumber\\
&=&
o(n^{-1/2}).
\end{eqnarray}
Finally, we have
\begin{eqnarray}\label{step7}
O\left(\frac{\log^9 n}{n^2}\right) \sum_{0\leq r\leq\log^{3/2} n}\frac{e^{rc_n}}{r!} &\leq&
O\left(\frac{\log^9 n}{n^2}\right)e^{c_n}\nonumber\\
&=&o(n^{-1/2}).
\end{eqnarray}
Together, (\ref{step1}), (\ref{step2}), (\ref{step3}), (\ref{step4}), (\ref{step5}), (\ref{step6})
and (\ref{step7}) prove (\ref{e1}).

To show (\ref{e2}), let $Y$ be the number of pairs $S_i$, $S_j$ in
an partition in $\cals_n$ chosen uniformly at random for which
$\psi(S_i)=\psi(S_j)$. For such $S_i$, $S_j$ of size $|S_i|=|S_j|=k$,
 the probability that
they are present in the random partition is $B(2n-2k)/B(2n)$.
The total number of pairs $S_i$, $S_j$ of size $k$ is bounded
by ${n\choose k} 2^k$ (the number of ways of choosing a subset $J$ of size
$k$ from $[n]$ times a bound on the number of ways of choosing two subsets
$S_1$, $S_2$ of $[2n]$ of size $k$
such that either $j\in S_1$ and $j+n\in S_2$
or $j+n\in S_1$ and $j\in S_2$ for all $j\in J$.)
Therefore, using (\ref{ratiobound}) we get
\begin{eqnarray*}
1-\frac{|E_{2,n}|}{B_{2n}}&=&
\BP(Y>0)\\
&\leq& \BE Y\\
&\leq &
\sum_{k=1}^n{n\choose k}2^k\frac{B_{2n-2k}}{B_{2n}}\\
&\leq&
\sum_{k=1}^n{n\choose k}2^k\frac{(C\log 2n)^{2k}}{(2n)_{2k}}\\
&\leq&
\sum_{k=1}^n\frac{(n)_k(2C^2\log^2 2n))^{k}}{(2n)_{2k}k!}\\
&=&
O\left(\frac{\log^2 n}{n}\right).
\end{eqnarray*}
\hfill\qed

Lemma~\ref{asymp} and (\ref{tnexp}) along with
$$
\frac{|C_n|}{B_{2n}}\leq \frac{|E_{1,n}|}{B_{2n}}
$$
and
$$
\frac{|C_n|}{B_{2n}}\geq \frac{|E_{1,n}|-(B_{2n}-|E_{2,n}|)}{B_{2n}}
$$
prove (\ref{tasymp}) for $t_n$.

\subsection{Asymptotic enumeration of 2-covers}
In this subsection we prove (\ref{tasymp}) for $s_n$. Recall that $U_n$ denotes the set of
2-covers of $[n]$. Each element of $E_{1,n}$ is mapped to a unique $\bfa\in U_n$ by $\phi$. Given
$\omega=\{\tilde{S}_1,\tilde{S}_2,\ldots,\tilde{S}_m\}\in\cals_n$, let $Z(\omega)$ be the number
of pairs $\{i_1,i_2\}$ such that $\psi(\tilde{S}_{i_1})=\psi(\tilde{S}_{i_2})$. Note that in the
case $\omega\in E_{1,n}$ we have $Z(\omega)=\rho$ with $\rho$ defined with respect to
$\bfa=\phi(\omega)$ in the statement of Lemma~\ref{improper}.

Define $D_{\rho,n}$ for $\rho\in\{0,1,\ldots,n\}$ to be
$$D_{\rho,n}=\{\omega\in E_{1,n}:Z(\omega)=\rho\}.$$
Note that $D_{0,n}=C_n$. By Lemma~\ref{improper},
\begin{eqnarray*}
u_n&=&\sum_{\rho=0}^n |D_{\rho,n}| 2^{-n+\rho}\\
&=&
|C_n|2^{-n}+ \sum_{\rho=1}^n |D_{\rho,n}|2^{\rho}\\
&=& B_{2n} 2^{-n} \left(\frac{|C_n|}{B_{2n}}+\sum_{\rho=1}^n \frac{|D_{\rho,n}|}{B_{2n}}
2^{\rho}\right).
\end{eqnarray*}
We have shown in the previous section that $C_n/B_{2n}\sim \sqrt{\log n/2n}$. Observe that
$\sum_{\rho=1}^n |D_{\rho,n}| 2^{\rho}/B_{2n}\leq\sum_{\rho=1}^n \BP(Z=\rho) 2^{\rho}$, where $Z$
was defined in the last paragraph and $\omega$ is chosen uniformly at random from ${\cals}_n$. In
light of these observations, to prove (\ref{tasymp}) for $s_n$ it suffices to prove that
\eq\label{tobeshown} \sum_{\rho=1}^n \BP(Z=\rho) 2^{\rho}=o\left(\sqrt{\frac{\log n}{2n}}\right).
\en

The quantity $\BP(Z\geq\rho)$ is equal to the probability that the randomly chosen element of
$\cals_n$ contains at least $\rho$ disjoint pairs of equal sets, therefore,
$$
\BP(Z\geq\rho)\leq
\sum_{s_1=1}^n\sum_{s_2=1}^n\cdots\sum_{s_\rho=1}^n
{n\choose {s_1,s_2,\ldots,s_\rho,n-\sum s_i}}
\frac{B_{2n-2\sum s_i}}{B_{2n}}
$$
Let $\sigma$ be defined by $\sigma=\sum_{i=1}^\rho s_i$. We can assume $\sigma\leq n$. From
(\ref{ratiobound}) we have
\begin{eqnarray*}
\BP(Z\geq\rho)&\leq&
\sum_{s_1=1}^n\sum_{s_2=1}^n\cdots\sum_{s_\rho=1}^n
{n\choose {s_1,s_2,\ldots,s_\rho,n-\sigma}}
\frac{(C\log n)^{2\sigma}}{(2n)_{2\sigma}}\\
&=&
\sum_{s_1=1}^n\sum_{s_2=1}^n\cdots\sum_{s_\rho=1}^n
\frac{(n)_\sigma}{\prod_i s_i!}
\frac{(C\log n)^{2\sigma}}{(2n)_{2\sigma}}.
\end{eqnarray*}
Observing that
$$
\frac{(n)_\sigma}{(2n)_{2\sigma}}=
\frac{(n)_\sigma}{(2n)_\sigma(2n-\sigma)_\sigma}\leq
\frac{1}{(2n)_\sigma}\leq n^{-\sigma},
$$
we have
\begin{eqnarray*}
\BP(Z\geq\rho)&\leq&
\sum_{\sigma=\rho}^n
\sum_{\stackrel{s_1,\ldots,s_\rho:}{\sum_i s_i=\sigma}}
\frac{1}{\prod_i s_i!}
\left(\frac{C^2\log^2 n}{n}\right)^\sigma\\
&=&
\sum_{\sigma=\rho}^n
\frac{\rho^\sigma}{\sigma!}
\left(\frac{C^2\log^2 n}{n}\right)^\sigma\\
\end{eqnarray*}
Therefore,
\begin{eqnarray*}
\sum_{\rho=1}^n \BP(Z=\rho) 2^{\rho}&\leq&
\sum_{\rho=1}^n \BP(Z\geq\rho) 2^{\rho}\\
&\leq&
\sum_{\rho=1}^n \sum_{\sigma=\rho}^n \frac{2^\rho \rho^\sigma}{\sigma!}
\left(\frac{C^2\log^2 n}{n}\right)^\sigma\\
&=&
\sum_{\sigma=1}^n \sum_{\rho=1}^\sigma \frac{2^\rho \rho^\sigma}{\sigma!}
\left(\frac{C^2\log^2 n}{n}\right)^\sigma\\
&\leq&
\sum_{\sigma=1}^n \sum_{\rho=1}^\sigma \frac{\rho^\sigma}{\sigma!}
\left(\frac{2C^2\log^2 n}{n}\right)^\sigma\\
&\leq&
\sum_{\sigma=1}^n \frac{(\sigma+1)^\sigma}{\sigma!}
\left(\frac{2C^2\log^2 n}{n}\right)^\sigma\\
&=&
O\left(
\frac{\log^2 n}{n}\right)\\
&=&
o\left(\sqrt{\frac{\log n}{2n}}\right).
\end{eqnarray*}
The last estimate proves (\ref{tobeshown}).
\hfill\qed

\section{Restricted 2-covers and line graphs: an analytic approach}

Our proof of (\ref{uasymp}) will use generating function analysis.
Let $a_{n,m}$ be the number of restricted, proper 2-covers on $[n]$ with $m$ blocks.
The generating function for restricted, proper $2$-covers
$$
A(x,y)=\sum_{n=0}^\infty \sum_{m=1}^{2n} \frac{a_{n,m}}{n!} x^n y^m
$$
equals \eq\label{genfunc} A(x,y)=\exp\left(-y-\frac{xy^2}{2}\right) \sum_{m\geq
0}\frac{y^m}{m!}(1+x)^{{m\choose 2}}; \en see page 203 of \cite{GJ}. Therefore, \eq\label{Vexp}
V(x)=A(x,1)= e^{-1}\sum_{m=0}^\infty\frac{1}{m!}(1+x)^{{m\choose 2}}e^{-x/2} \en and
\eq\label{vsum} v_n= n!e^{-1} \sum_{m=0}^\infty\frac{m^{2n}}{m!}
\sum_{k=0}^n\frac1{k!}\left(-\frac{1}2\right)^k m^{-2n}{{m\choose 2}\choose n-k}. \en Note that
for $m\geq 2$,
\begin{eqnarray}\label{sumbound}
\left|\sum_{k=0}^n\frac{n!}{k!}\left(-\frac{1}{2}\right)^k m^{-2n}{{m\choose 2}\choose n-k}\right|
&\leq& \sum_{k=0}^n{n\choose k}\left(\frac{1}{2}\right)^km^{-2n}{m\choose 2}^{n-k}
\nonumber\\
&\leq&
2^{-n}\sum_{k=0}^n{n\choose k}m^{-2k}\nonumber\\
&\leq& 2^{-n}\left(\frac{1+m^{-2}}{2}\right)^n=O(2^{-n}).
\end{eqnarray}

We will make use of the asymptotic analysis of the Bell numbers
in Example~5.4 of \cite{O}, which uses the identity
$$
B_n=e^{-1}\sum_{m=0}^\infty\frac{m^n}{m!}.
$$
Let $m_0$ be the nearest integer to $\frac{2n}{W(2n)}$, where $W$ is defined by
(\ref{Lambertdef}). (The choice of $m_0$ is slightly different here than in \cite{O}, but the
analysis giving (\ref{outer}) and (\ref{inner}) below remains valid.) In \cite{O} it is proved
that \eq\label{outer} \sum_{\stackrel{1\leq m\leq n}{|m-m_0|>\sqrt{n}\log n}} \frac{m^{2n}}{m!}=
O\left(\frac{m_0^{2n}}{m_0!}\sqrt{n}\exp\left(-(\log n)^3\right)\right) \en and that
\begin{eqnarray}\label{inner} \sum_{\stackrel{1\leq m\leq n}{|m-m_0|\leq\sqrt{n}\log n}}
\frac{m^{2n}}{m!}&=&
\frac{m_0^{2n+1}}{m_0!}\sqrt{\frac{2\pi}{2n+m_0}} \left(1+O\left((\log n)^6
n^{-1/2}\right)\right)\\
&\sim&eB_{2n}. \end{eqnarray}

It follows from (\ref{sumbound}) and (\ref{outer}) that
\begin{eqnarray}
\sum_{\stackrel{1\leq m\leq n}{|m-m_0|>\sqrt{n}\log n}} \frac{m^{2n}}{m!}
\sum_{k=0}^n\frac{n!}{k!}\left(-\frac{1}2\right)^k m^{-2n}{{m\choose 2}\choose n-k}
&=&O\left(\frac{m_0^{2n}}{m_0!}\sqrt{n}2^{-n}\exp\left(-(\log n)^3\right)\right)\nonumber\\
&=&O\left(B_{2n}2^{-n}\exp\left(-\frac{(\log n)^3}{2}\right)\right)\label{outer2}. \end{eqnarray}
We have \eq\label{decomp} \sum_{\stackrel{1\leq m\leq n}{|m-m_0|\leq\sqrt{n}\log n}}
\frac{m^{2n}}{m!} \sum_{k=0}^n\frac{n!}{k!}\left(-\frac{1}2\right)^k m^{-2n}{{m\choose 2}\choose
n-k} = \sum_{\stackrel{1\leq m\leq n}{|m-m_0|\leq\sqrt{n}\log n}} \frac{m^{2n}}{m!}
m^{-2n}n!{{m\choose 2}\choose n}+\Delta, \en where
$$
\Delta:=
\sum_{\stackrel{1\leq m\leq n}{|m-m_0|\leq\sqrt{n}\log n}}
\frac{m^{2n}}{m!}
\sum_{k=1}^n\frac{n!}{k!}\left(-\frac{1}2\right)^k
m^{-2n}{{m\choose 2}\choose n-k}
$$
is bounded by
\begin{eqnarray*}
|\Delta|&\leq&
\sum_{\stackrel{1\leq m\leq n}{|m-m_0|\leq\sqrt{n}\log n}}
\frac{m^{2n}}{m!}
\sum_{k=1}^n\frac{n!}{k!}
m^{-2n}{{m\choose 2}\choose n}\left(\frac{n}{{m\choose 2}-n}\right)^k\\
&=&
O\left(\frac{\log^2 n}{n}\right)
\sum_{\stackrel{1\leq m\leq n}{|m-m_0|\leq\sqrt{n}\log n}}
\frac{m^{2n}}{m!}
m^{-2n} n!{{m\choose 2}\choose n}.
\end{eqnarray*}
One may show that uniformly for $m$ in the range $|m-m_0|\leq\sqrt{n}\log n$ $$ m^{-2n}{{m\choose
2}\choose n} n!= 2^{-n}\exp\left(-{n\over m_0}-{n^2\over m_0^2}\right)
\left(1+O\left(n^{-1/2}\log^6 n\right)\right), $$ hence, \eq\label{Deltabound}
|\Delta|=O\left(\frac{\log^2 n}{n}\right) 2^{-n}\exp\left(-{n\over m_0}-{n^2\over m_0^2}\right)
B_{2n}. \en The main term of (\ref{decomp}) is
\begin{eqnarray}\label{mainasymp}
\sum_{\stackrel{1\leq m\leq n}{|m-m_0|\leq\sqrt{n}\log n}} \frac{m^{2n}}{m!} m^{-2n}n!{{m\choose
2}\choose n} &=& 2^{-n}\exp\left(-{n\over m_0}-{n^2\over m_0^2}\right) \left(1+o(1)\right)
\sum_{\stackrel{1\leq m\leq n}{|m-m_0|\leq\sqrt{n}\log n}}
\frac{m^{2n}}{m!}\nonumber\\
&=& eB_{2n}2^{-n}\exp\left(-{n\over m_0}-{n^2\over m_0^2}\right)
(1+o(1))\nonumber\\
&=& eB_{2n}\frac1{2^n\sqrt n}e^{-\left(\frac12\log(2n/\log n)\right)^2} (1+o(1))
\end{eqnarray}
where we have used the asymptotic expansion (\ref{Lambertasymp}) and the definition of $m_0$ at
the last step. Now (\ref{vsum}), (\ref{outer2}), (\ref{Deltabound}) and (\ref{mainasymp}) prove
(\ref{uasymp}) for $v_n$.

In the previous argument the result would have been the same
if the $e^{-x/2}$
in (\ref{Vexp}) were replaced by $1$
because in the Taylor expansion
of $e^{-x/2}$ the constant term $1$ corresponds to
the main term of (\ref{decomp})
and the higher order terms contribute to $\Delta$, which is negligable.
The argument for restricted partitions and line graphs
are similar, starting
from the identities obtained from Proposition~\ref{genfunc} and
(\ref{Vexp})
$$
U(x)=
e^{-1}\sum_{m=0}^\infty\frac{1}{m!}(1+x)^{{m\choose 2}}e^{x/2}.
$$
and
$$
L(x)=
e^{-1}\sum_{m=0}^\infty\frac{1}{m!}(1+x)^{{m\choose 2}}e^{x/2-x^3/6}.
$$
In each case only the contribution of the constant term of the Taylor
expansion of the exponential is $1$ and the remaining terms contribute
to a quantity like $\Delta$ which is asymptotically insignificant.
\hfill\qed

\end{document}